\documentclass[12pts]{amsart}
\usepackage{amsmath, amssymb, graphicx}
\usepackage{amsthm}
\usepackage[arrow,matrix]{xy}
\begin{document}
\def\eq#1{{\rm(\ref{#1})}}
\theoremstyle{plain}
\newtheorem{thm}{Theorem}[section]
\newtheorem{lem}[thm]{Lemma}
\newtheorem{prop}[thm]{Proposition}
\newtheorem{cor}[thm]{Corollary}
\theoremstyle{definition}
\newtheorem{dfn}[thm]{Definition}
\newtheorem{rem}[thm]{Remark}
\def\Ker{\mathop{\rm Ker}}
\def\Coker{\mathop{\rm Coker}}
\def\ind{\mathop{\rm ind}}
\def\Re{\mathop{\rm Re}}
\def\vol{\mathop{\rm vol}}
\def\SO{\mathbin{\rm SO}}
\def\Im{\mathop{\rm Im}}
\def\min{\mathop{\rm min}}
\def\Spec{\mathop{\rm Spec}\nolimits}
\def\Hol{{\textstyle\mathop{\rm Hol}}}
\def\ge{\geqslant}
\def\le{\leqslant}
\def\C{{\mathbin{\mathbb C}}}
\def\R{{\mathbin{\mathbb R}}}
\def\N{{\mathbin{\mathbb N}}}
\def\Z{{\mathbin{\mathbb Z}}}
\def\D{{\mathbin{\mathcal D}}}
\def\H{{\mathbin{\mathcal H}}}
\def\M{{\mathbin{\mathcal M}}}
\def\al{\alpha}
\def\be{\beta}
\def\ga{\gamma}
\def\de{\delta}
\def\ep{\epsilon}
\def\io{\iota}
\def\ka{\kappa}
\def\la{\lambda}
\def\ze{\zeta}
\def\th{\theta}
\def\vp{\varphi}
\def\si{\sigma}
\def\up{\upsilon}
\def\Up{\Upsilon}
\def\vars{\varsigma}
\def\varr{\varrho}
\def\om{\omega}
\def\De{\Delta}
\def\Ga{\Gamma}
\def\Th{\Theta}
\def\La{\Lambda}
\def\Om{\Omega}
\def\ts{\textstyle}
\def\sst{\scriptscriptstyle}
\def\sm{\setminus}
\def\na{\nabla}
\def\pd{\partial}
\def\op{\oplus}
\def\ot{\otimes}
\def\bigop{\bigoplus}
\def\iy{\infty}
\def\ra{\rightarrow}
\def\longra{\longrightarrow}
\def\dashra{\dashrightarrow}
\def\t{\times}
\def\w{\wedge}
\def\d{{\rm d}}
\def\bs{\boldsymbol}
\def\ci{\circ}
\def\ti{\tilde}
\def\ov{\overline}
\def\md#1{\vert #1 \vert}
\def\nm#1{\Vert #1 \Vert}
\def\bmd#1{\big\vert #1 \big\vert}
\def\cnm#1#2{\Vert #1 \Vert_{C^{#2}}} 
\def\lnm#1#2{\Vert #1 \Vert_{L^{#2}}} 
\def\bnm#1{\bigl\Vert #1 \bigr\Vert}
\def\bcnm#1#2{\bigl\Vert #1 \bigr\Vert_{C^{#2}}} 
\def\blnm#1#2{\bigl\Vert #1 \bigr\Vert_{L^{#2}}} 

\title[DEFORMATIONS OF AC COASSOCIATIVE 4-FOLDS]{Deformations of
Asymptotically Cylindrical Coassociative Submanifolds with Moving
Boundary}
\author{Sema Salur}
\address {Department of Mathematics, University of Rochester, Rochester, NY, 14627.}
\email{salur@math.rochester.edu }

\begin{abstract} In an earlier paper, \cite{Salu0}, we proved that
given an asymptotically cylindrical $G_2$-manifold $M$ with a
Calabi--Yau boundary $X$, the moduli space of coassociative
deformations of an asymptotically cylindrical coassociative 4-fold
$C\subset M$ with a fixed special Lagrangian boundary $L\subset X$
is a smooth manifold of dimension $\dim V_+$, where $V_+$ is the
positive subspace of the image of $H^2_{\rm cs}(C,\mathbb{R})$ in
$H^2(C,\mathbb{R})$. In order to prove this we used the powerful
tools of Fredholm Theory for noncompact manifolds which was
developed by Lockhart and McOwen \cite{LM1}, \cite{LM2}, and
independently by Melrose \cite{Mel1}, \cite{Mel2}.

In this paper, we extend our result to the moving boundary case.
Let $\Up:H^2(L,\mathbb{R})\rightarrow H^3_{\rm cs}(C,\mathbb{R})$
be the natural projection, so that $\ker\Up$ is a vector subspace
of $H^2(L,\mathbb{R})$. Let $F$ be a small open neighbourhood of 0
in $\ker\Up$ and $L_s$ denote the special Lagrangian submanifolds
of $X$ near $L$ for some $s\in F$ and with phase $i$. Here we
prove that the moduli space of coassociative deformations of an
asymptotically cylindrical coassociative submanifold $C$
asymptotic to $L_s\times (R,\infty)$, $s\in F$, is a smooth
manifold of dimension equal to $\dim V_++\dim\ker\Up=\dim V_+
+b^2(L)-b^0(L)+b^3(C)-b^1(C)+b^0(C)$.


\end{abstract}

\maketitle
\section{Introduction}

Let $(M,\vp,g_M)$ be a connected, complete, asymptotically
cylindrical $G_2$-manifold with a $G_2$-structure $(\vp,g_M)$ and
asymptotic to $X\t(R,\iy)$, $R>0$, with decay rate $\al<0$, where
$X$ is a Calabi--Yau $3$-fold. An asymptotically cylindrical
$G_2$-manifold $M$ is a noncompact Riemannian 7-manifold with zero
Ricci curvature whose holonomy group Hol$(g_M)$ is a subgroup of
the exceptional Lie group $G_2$. It is equipped with a covariant
constant 3-form $\vp$ and a 4-form $*\vp$. There are two natural
classes of noncompact calibrated submanifolds inside $M$
corresponding to $\vp$ and $*\vp$ which are called asymptotically
cylindrical associative 3-folds and coassociative 4-folds,
respectively.

\vspace{.1in}

The Floer homology programs for asymptotically cylindrical
Calabi-Yau and $G_2$-manifolds lead to the construction of brand
new Topological Quantum Field Theories. In order to construct
consistent TQFT's, the first step is to understand the
deformations of asymptotically cylindrical calibrated submanifolds
with some small decay rate inside Calabi-Yau and $G_2$-manifolds.
A fundamental question is whether these noncompact submanifolds
have smooth deformation spaces.

\vspace{.1in}

For this purpose, in \cite{Salu0} we studied the deformation space
of an asymptotically cylindrical coassociative submanifold $C$ in
an asymptotically cylindrical $G_2$-manifold $M$ with a Calabi-Yau
boundary $X$. We assumed that the boundary $\partial C =L$ is a
special Lagrangian submanifold of $X$. Using the analytic set-up
which was developed for elliptic operators on asymptotically
cylindrical manifolds by Lockhart-McOwen and Melrose, \cite{LM1},
\cite{LM2}, \cite{Mel1}, \cite{Mel2}, we proved that for fixed
boundary $\partial C = L$ this moduli space is a smooth manifold.

\vspace{.1in}

In order to develope the Floer homology program for coassociative
submanifolds the next step is to parametrize coassociative
deformations with moving boundary. In this paper, we extend our
previous result in \cite{Salu0} to the moving boundary case and
show that if the special Lagrangian boundary is allowed to move,
one still obtains a smooth moduli space. Moreover, we can
determine the dimension of this space. These two technical results
on fixed and moving boundary cases then parametrize all
asymptotically cylindrical coassociative deformations. Remarkably,
using this parametrization of deformations and basic algebraic
topology we can verify one of the main claims of Leung in \cite
{Leun2} to prove that the boundary map from the moduli space of
coasssociative cycles to the moduli space of special Lagrangian
cycles is a Lagrangian immersion. With this result in hand, one
can start defining the Floer homology of coassociative cycles with
special Lagrangian boundary, which then leads to the construction
of the TQFT of these cycles.

\vspace{.1in}

 In this paper we prove the following theorem.

\begin{thm}
Let $(M,\vp,g)$ be an asymptotically cylindrical $G_2$-manifold
asymptotic to $X\t(R,\iy)$, $R>0$, with decay rate $\al<0$, where
$X$ is a Calabi--Yau $3$-fold with metric $g_X$. Let $C$ be an
asymptotically cylindrical coassociative $4$-fold in $M$
asymptotic to $L\t(R',\iy)$ for $R'>R$ with decay rate $\be$ for
$\al\le\be<0$, where $L$ is a special Lagrangian $3$-fold in $X$
with phase $i$, and metric $g_L=g_X\vert_L$.

Let also $\Up:H^2(L,\mathbb{R})\rightarrow H^3_{cs}(C,\mathbb{R})$
be the natural projection, so that $\ker\Up$ is a vector subspace
of $H^2(L,\mathbb{R})$ and let $F$ be a small open neighbourhood
of 0 in $\ker\Up$. Let also $L_s$ be the special Lagrangian
submanifolds of $X$ near $L$ for some $s\in F$ and with phase $i$.
Then for any sufficiently small $\ga$ the moduli space $\M_C^\ga$
of asymptotically cylindrical coassociative submanifolds in $M$
close to $C$, and asymptotic to $L_s\t(R',\iy)$ with decay rate
$\ga$, is a smooth manifold of dimension $\dim
V_++\dim\ker\Up=\dim V_+ +b_2(L)-b_0(L)+b_3(C)-b_1(C)+b_0(C)$,
where $V_+$ is the positive subspace of the image of $H^2_{\rm
cs}(C,\R)$ in~$H^2(C,\R)$. \label{co1thm1}
\end{thm}

\begin{rem}
In Theorem \ref{co1thm1} above, as in the fixed boundary case,
\cite{Salu0}, we require $\ga$ to satisfy $\be<\ga<0$, and
$(0,\ga^2]$ to contain no eigenvalues of the Laplacian $\De_L$ on
functions on $L$, and $[\ga,0)$ to contain no eigenvalues of the
operator $-*\d$ on coexact $1$-forms on $L$. These hold provided
$\ga<0$ is small enough. This assumption is needed to guarantee
that the linearized operator of the deformation map for
asymptotically cylindrical coassociative submanifolds is Fredholm.
\end{rem}

\begin{rem}
In our previous paper, \cite{Salu0}, we proved that the dimension
of the deformation space of asymptotically cylindrical
coassociative submanifold $C$ with fixed special Lagrangian
boundary $L$ is given as $V_+$ where $V_+$ is the positive
subspace of the image of $H^2_{\rm cs}(C,\R)$ in $H^2(C,\R)$. One
should note that in Theorem \ref{co1thm1} above, in the case when
$L$ is a special Lagrangian homology 3-sphere, $b_1(L)=0$ and
hence $L$ is rigid and there are no special Lagrangian
deformations of $L$ so it behaves like a fixed boundary. Then one
can use basic algebraic topology to show that for this special
case Theorem \ref{co1thm1} gives us $\dim V_++\dim\ker\Up=\dim V_+
+b_2(L)-b_0(L)+b_3(C)-b_1(C)+b_0(C)=\dim V_+$ which is consistent
with the result of our previous paper, \cite{Salu0}, for fixed
boundary.

\end{rem}

\begin{rem}
Theorem \ref{co1thm1} applies to examples of asymptotically
cylindrical $G_2$-manifolds constructed by Kovalev, in
\cite{Kova}. These 7-manifolds are of the form $X\times S^1$ and
are obtained by the direct product of asymptotically cylindrical
6-manifolds $X$ with holonomy $SU(3)\subset G_2$ and the circle
$S^1$. It is still an open question whether there exist
asymptotically cylindrical 7-manifolds with holonomy group $G_2$.
\end{rem}

\vspace{.1in}


The outline of the paper is as follows: In \S\ref{yco1}, we begin
with a general discussion of $G_2$ geometry, asymptotically
cylindrical manifolds and basic tools used in elliptic theory such
as weighted Sobolev spaces. In \S\ref{yco2}, we give an overview
of deformations of an asymptotically cylindrical coassociative
submanifold with fixed special Lagrangian boundary, followed by a
sketch proof of our main theorem in \cite{Salu0}, where we showed
that the moduli space of such deformations is smooth and
calculated its dimension. In \S\ref{yco3} we introduce the
analytic set-up for deformations with moving (free) boundary and
prove Theorem \ref{co1thm1}. Finally, in \S\ref{yco4} using
Theorem \ref{co1thm1}, we verify one of the main claims of Leung
in \cite{Leun2}, which is necessary to prove that the boundary map
from the moduli space of coassociative cycles into the moduli
space of special Lagrangian cycles is a Lagrangian immersion.

\section{$\mathbf{G_2}$-manifolds and coassociative submanifolds}

\label{yco1}

We now explain some background material on asymptotically
cylindrical $G_2$-manifolds and their coassociative submanifolds.
We will also review the Fredholm Theory for elliptic operators and
weighted Sobolev spaces. A good reference for $G_2$ geometry is
Harvey and Lawson \cite{HaLa}. The details of the elliptic theory
on noncompact manifolds can be found in
\cite{LM1},\cite{LM2},\cite{Mel1},\cite{Mel2}.

\vspace{.1in}

\subsection{$\mathbf{G_2}$ Geometry}

The imaginary octonions $Im \mathbb{O} =\mathbb {R}^7$ is equipped
with the cross product $\times: \mathbb {R}^7\times \mathbb {R}^7
\to \mathbb {R}^7$ defined by $u\times v=Im(u\cdot\bar{v})$ where
$\cdot$ is the octonionic multiplication. The exceptional Lie
group $G_{2}$ can be defined as the linear automorphisms of
Im$\mathbb{O}$ preserving this cross product operation. It is a
compact, semisimple, and 14-dimensional subgroup of $SO(7)$.

\begin{dfn}
A smooth $7$-manifold $M$ has a $G_2$-structure if its tangent
frame bundle reduces to a $G_2$ bundle.
\end{dfn}

It is known that if $M$ has a $G_2$-structure then there is a
$G_2$-invariant 3-form $\varphi \in \Omega^{3}(M)$ which can be
written in an orthonormal frame as

$$\varphi=e^1 \wedge e^2 \wedge e^3+e^1 \wedge e^4 \wedge e^5+e^1\wedge e^6 \wedge e^7+e^2 \wedge e^4 \wedge e^6-e^2\wedge e^5
\wedge e^7$$
$$-e^3 \wedge e^4 \wedge e^7-e^3 \wedge e^5 \wedge e^6. $$

\vspace{.1in}

This $G_2$-invariant 3-form $\varphi$ gives an orientation $\mu
\in \Omega^{7}(M)$ on $M$ and $\mu$ determines a metric
$g=g_{\varphi }= \langle \;,\;\rangle$ on $M$ defined as
\begin{equation*}
\langle u,v \rangle=[ i_{u}(\varphi ) \wedge i_{v}(\varphi )\wedge
\varphi  ]/\mu ,
\end{equation*}
\noindent where $i_{v}=v\lrcorner $ is the interior product with a
vector $v$. So from now on we will refer to the term $(\varphi,
g)$ as a $G_2$-structure and the term $(M,\varphi, g)$ as a
manifold with $G_2$-structure.

\vspace{.1in}

\begin{dfn}

$(M,\varphi, g)$ is a $G_2$-manifold if  $\nabla\varphi=0$, i.e.
the $G_2$-structure $\varphi$ is torsion-free.
\end{dfn}

One can show that for a manifold with $G_2$-structure
$(M,\varphi,g)$, the following are equivalent:

\begin{itemize}
\setlength{\parsep}{0pt} \setlength{\itemsep}{0pt} \item[{\rm(i)}]
$\Hol(g)\subseteq G_2$, \item[{\rm(ii)}] $\na\vp=0$ on $M$, where
$\na$ is the Levi-Civita connection of $g$, \item[{\rm(iii)}]
$\d\vp=0$ and $\d^*\vp=0$ on $M$.
\end{itemize}

\vspace{.1in}

Harvey and Lawson, \cite{HaLa}, showed that there are minimal
submanifolds of $G_2$-manifolds calibrated by $\varphi$ and
$*\varphi$.

\begin{dfn}
Let $(M, \varphi,g )$ be a $G_2$-manifold. A 4-dimensional
submanifold $C\subset M$ is called {\em coassociative } if
$\varphi|_C=0$. A 3-dimensional submanifold $Y\subset M$ is called
{\em associative} if $\varphi|_Y\equiv vol(Y)$; this condition is
equivalent to $\chi|_Y\equiv 0$,  where $\chi \in \Omega^{3}(M,
TM)$ is the  tangent bundle valued 3-form defined by the identity:

\begin{equation*}
\langle \chi (u,v,w) , z \rangle=*\varphi  (u,v,w,z)
\end{equation*}

\end{dfn}

\vspace{.05in}
\subsection{Asymptotically Cylindrical $\mathbf{G_2}$-Manifolds and Coassociative Submanifolds}
\label{subsec2}

 Next, we recall basic properties of asymptotically
cylindrical $G_2$-manifolds and the coassociative 4-folds. We also
need these definitions and the analytic set-up in Section
\ref{yco3}. More details on the subject can be found in
\cite{Salu0}.

\vspace{.1in}

\begin{dfn}
A $G_2$-manifold $(M_0,\vp_0,g_0)$ is called cylindrical if
$M_0=X\t\R$ and $(\vp_0,g_0)$ is compatible with the product
structure, that is,

\begin{equation*}
\vp_0=\Re\Om+\om\w\d t \qquad\text{and}\qquad g_0=g_X+\d t^2,
\end{equation*}

\noindent where $X$ is a connected, compact Calabi--Yau 3-fold
with K\"{a}hler form $\om$, Riemannian metric $g_X$ and
holomorphic (3,0)-form $\Om$.

\end{dfn}

\vspace{.1in}

\begin{dfn}
Let $X$ be a Calabi--Yau $3$-fold with metric $g_X$. $(M,\vp,g)$
is called an asymptotically cylindrical $G_2$-manifold asymptotic
to $X\t(R,\iy)$, $R>0$, with decay rate $\al<0$ if there exists a
cylindrical $G_2$-manifold $(M_0,\vp_0,g_0)$ with $M_0=X\t\R$, a
compact subset $K\subset M$, a real number $R$, and a
diffeomorphism $\Psi:X\t(R,\iy)\ra M\sm K$ such that
$\Psi^*(\vp)=\vp_0+\d\xi$ for some smooth 2-form $\xi$ on
$X\t(R,\iy)$ with $\bmd{\na^k\xi}=O(e^{\al t})$ on $X\t(R,\iy)$
for all $k\ge 0$, where $\na$ is the Levi-Civita connection of the
cylindrical metric $g_0$.

\end{dfn}

\vspace{.1in}

An asymptotically cylindrical $G_2$-manifold $M$ has one end
modelled on $X\t(R,\iy)$, and as $t\ra\iy$ in $(R,\iy)$ the
$G_2$-structure $(\vp,g)$ on $M$ converges to order $O(e^{\al t})$
to the cylindrical $G_2$-structure on $X\t(R,\iy)$, with all of
its derivatives. As in \cite{Salu0}, we suppose $M$ and $X$ are
connected, that is, we allow $M$ to have only one end. We showed
earlier that an asymptotically cylindrical $G_2$-manifold can have
at most one cylindrical end, or otherwise the holonomy reduces,
\cite{Salu01}.

\vspace{.1in}

We can also define calibrated submanifolds of noncompact
$G_2$-manifolds. In \cite{Salu0}, we introduced definitions of
cylindrical and asymptotically cylindrical coassociative
submanifolds of a $G_2$-manifold.

\vspace{.1in}

\begin{dfn}
Let $(M_0,\vp_0,g_0)$ be a cylindrical $G_2$-manifold. A 4-fold
$C_0$ is called cylindrical coassociative submanifold of $M_0$ if
$C_0=L\times\mathbb{R}$ for some compact special Lagrangian 3-fold
$L$ with phase $i$ in the Calabi-Yau manifold $X$.
\end{dfn}

\begin{dfn}
Let $M$ be an asymptotically cylindrical $G_2$-manifold asymptotic
to $X\t(R,\iy)$, $R>0$, with decay rate $\al<0$.  Let also $C$ be
a connected, complete asymptotically cylindrical $4$-fold in $M$
asymptotic to $L\t(R',\iy)$ for $R'>R$ with decay rate $\be$ for
$\al\le\be<0$ and let $L$ be a compact special Lagrangian $3$-fold
in $X$ with phase $i$, and metric $g_L=g_X\vert_L$.  Then $C$ is
called asymptotically cylindrical coassociative 4-fold if there
exists a compact subset $K'\subset C$, a normal vector field $v$
on $L\t(R',\iy)$ for some $R'>R$, and a diffeomorphism
$\Phi:L\t(R',\iy)\ra C\sm K'$ such that the following diagram
commutes:

\begin{equation}
\begin{gathered}
\xymatrix{X\t (R',\iy) \ar[d]^\subset & L\t (R',\iy)
\ar[l]^{\exp_v} \ar[r]_{\Phi} & (C\sm K')
\ar[d]^\subset \\
X\t (R,\iy)\ar[rr]^\Psi && (M\sm K), }
\end{gathered}
\label{co2eq3}
\end{equation}

\vspace{.05in}

and $\bmd{\na^kv}=O(e^{\be t})$ on $L\t(R',\iy)$ for all $k\ge 0$.

\label{defn1}
\end{dfn}

\vspace{.1in}

Diagram (\ref{co2eq3}) implies that $C$ in $M$ is asymptotic to
the cylinder $C_0$ in $M_0=X\t\R$ as $t\ra\iy$ in $\R$, to order
$O(e^{\be t})$. In other words, $C$ can be written near infinity
as the graph of a normal vector field $v$ to $C_0=L\t\R$ in
$M_0=X\t\R$, so that $v$ and its derivatives are $O(e^{\be t})$.
Here we require $C$ but not $L$ to be connected, so it is possible
that $C$ to have multiple ends.

\vspace{.1in}

\subsection{Elliptic operators on noncompact manifolds} \label{subsec3} We will conclude this section with a brief review
of the weighted Sobolev spaces and the results of Lockhart and
McOwen, \cite{LM1}, \cite{LM2}, about Fredholm properties of
elliptic operators on manifolds with cylindrical ends.

\vspace{.1in}

Let $C$ and $L$ be as in Definition \ref{defn1} above. $E_0$ is
called a cylindrical vector bundle on $L\t\R$ if it is invariant
under translations in $\R$. Let $h_0$ be a smooth family of
metrics on the fibres of $E_0$ and $\na_{\sst E_0}$ a connection
on $E_0$ preserving $h_0$, with $h_0,\na_{\sst E_0}$ invariant
under translations in $\R$. Let $E$ be a vector bundle on $C$
equipped with metrics $h$ on the fibres, and a connection
$\na_{\sst E}$ on $E$ preserving $h$. We say that $E,h,\na_{\sst
E}$ are asymptotic to $E_0,h_0,\na_{\sst E_0}$ if there exists an
identification $\Phi_*(E)\cong E_0$ on $L\t(R',\iy)$ such that
$\Phi_*(h)=h_0+O(e^{\be t})$ and $\Phi_*( \na_{\sst E})=\na_{\sst
E_0}+O(e^{\be t})$ as $t\ra\iy$. Then we call $E,h,\na_{\sst E}$
asymptotically cylindrical.

\vspace{.1in}

\begin{dfn}
Let $\rho:C\ra\R$ be a smooth function satisfying $\Phi^*(\rho)
\equiv t$ on $L\t(R',\iy)$. For $p\geq 1$, $k\geq 0$ and
$\ga\in\R$ the weighted Sobolev space $L^p_{k,\ga}(E)$ is defined
to be the set of sections $s$ of $E$ that are locally integrable
and $k$ times weakly differentiable and for which the norm

\begin{equation}
\nm{s}_{L^p_{k,\ga}}=\Bigl(\sum_{j=0}^{k}\int_C
e^{-\ga\rho}\bmd{\na_{\sst E}^js}^p\d V\Bigr)^{1/p} \label{co3eq3}
\end{equation}

\noindent is finite.

\label{co3def1}
\end{dfn}

Note that the weighted Sobolev space, $L^p_{k,\ga}(E)$, is a
Banach space.

\vspace{.1in}

Now suppose $E,F$ are two asymptotically cylindrical vector
bundles on $C$, asymptotic to cylindrical vector bundles $E_0,F_0$
on $L\t\R$.

 \vspace{.1in}

\begin{dfn} Let $A_0:C^\iy(E_0)\ra C^\iy(F_0)$ be a cylindrical
elliptic operator of order $k$ invariant under translations in
$\mathbb R$. Let also $A:C^\iy(E)\ra C^\iy(F)$ be an elliptic
operator of order $k$ on $C$. $A$ is asymptotic to $A_0$ if under
the identifications $\Phi_*(E)\cong E_0$, $\Phi_*(F)\cong F_0$ on
$L\t(R',\iy)$ we have $\Phi_*(A)=A_0+O(e^{\be t})$ as $t\ra\iy$
for $\be<0$. Then $A$ is called an asymptotically cylindrical
elliptic operator.
\end{dfn}

\vspace{.1in}

It is well known that if $A$ is an elliptic operator on a compact
manifold then it should be Fredholm. But this is not the case for
the noncompact manifolds; there are examples of elliptic operators
which are not Fredholm, \cite{LM1}, \cite{LM2} and it turns out
that $A$ is Fredholm if and only if $\ga$ does not lie in a
discrete set $\D_{A_0}\subset\R$ which can be defined as follows:

\vspace{.1in}

\begin{dfn} Let $A$ and $A_0$ be elliptic operators on $C$ and $L\t\R$, so that $E,F$ have the same
fibre dimensions. Extend $A_0$ to the complexifications
$A_0:C^\iy(E_0\ot_\R\C)\ra C^\iy(F_0\ot_\R\C)$. Define $\D_{A_0}$
to be the set of $\ga\in\R$ such that for some $\de\in\R$ there
exists a nonzero section $s\in C^\iy(E_0\ot_\R\C)$ invariant under
translations in $\R$ such that $A_0(e^{(\ga+i\de)t}s)=0$.
\label{co3def3}
\end{dfn}

\vspace{.1in}

An important Fredholm property of elliptic operators on noncompact
manifolds has been shown by Lockhart and McOwen in
\cite[Th.~1.1]{LM2}:

\begin{thm} Let $(C,g)$ be an asymptotically cylindrical
Riemannian manifold asymptotic to $(L\t\R,g_0)$, and
$A:C^\iy(E)\ra C^\iy(F)$ an asymptotically cylindrical elliptic
operator on $C$ of order $k$ between asymptotically cylindrical
vector bundles $E,F$ on $C$, asymptotic to the cylindrical
elliptic operator $A_0: C^\iy(E_0)\ra C^\iy(F_0)$ on $L\t\R$. Let
$\D_{A_0}$ be defined as above.

Then $\D_{A_0}$ is a discrete subset of $\R$, and for $p>1$, $l\ge
0$ and $\ga\in\R$, the extension $A^p_{k+l,\ga}:L^p_{k+l,\ga}
(E)\ra L^p_{l,\ga}(F)$ is Fredholm if and only if
$\ga\notin\D_{A_0}$. \label{co3thm1}
\end{thm}

\vspace{.1in}

Note that Theorem \ref{co3thm1} plays an important role in our
choice of weighted Sobolev spaces in Theorem \ref{co1thm1}.

\section{Coassociative Deformations with Fixed Boundary}
\label{yco2}

Using the set-up in Section \ref{yco1}, we proved the following
theorem in \cite{Salu0}, for asymptotically cylindrical
coassociative deformations with fixed boundary.

\begin{thm} \cite[Thm 1.1.]{Salu0} Let $(M,\vp,g)$ be an asymptotically cylindrical $G_2$-manifold
asymptotic to $X\t(R,\iy)$, $R>0$, with decay rate $\al<0$, where
$X$ is a Calabi--Yau $3$-fold with metric $g_X$. Let $C$ be an
asymptotically cylindrical coassociative $4$-fold in $M$
asymptotic to $L\t(R',\iy)$ for $R'>R$ with decay rate $\be$ for
$\al\le\be<0$, where $L$ is a special Lagrangian $3$-fold in $X$
with phase $i$, and metric $g_L=g_X\vert_L$.

Then for some small $\ga$ the moduli space $\M_C^\ga$ of
asymptotically cylindrical coassociative submanifolds in $M$ close
to $C$, and asymptotic to $L\t(R',\iy)$ with decay rate $\ga$, is
a smooth manifold of dimension $\dim V_+$, where $V_+$ is the
positive subspace of the image of $H^2_{\rm cs}(C,\R)$ in
$H^2(C,\R)$. \label{co1thm}
\end{thm}

\begin{rem}

McLean proved the compact version of Theorem \ref{co1thm} in
\cite{McLe} and showed that the moduli space $\M_C$ of
coassociative $4$-folds isotopic to a compact coassociative
$4$-fold $C$ in $M$ is a smooth manifold of dimension $b^2_+(C)$.
There he modelled $\M_C$ on $\ti P^{-1}(0)$ for a nonlinear map
$\ti P$ between Banach spaces, whose linearization $\d\ti P(0,0)$
at 0 was the Fredholm map between Sobolev spaces

\begin{equation}
\d_++\d^*:L^p_{l+2}(\La^2_+T^*C)\t L^p_{l+2}(\La^4T^*C)\longra
L^p_{l+1}(\La^3T^*C). \label{co3eq1}
\end{equation}

McLean showed that $\d\ti P$ is onto the image of $\ti P$ and used
the Implicit Mapping Theorem for Banach spaces, \cite[Thm
1.2.5]{Joyc1} to conclude that $\ti P^{-1}(0)$ is smooth,
finite-dimensional and locally isomorphic to $\Ker\bigl((\d_++\d^*
)^p_{l+2}\bigr)$.

\end{rem}

\begin{proof}[Sketch proof of Theorem \ref{co1thm}]
Let $(M,\vp,g)$ be an asymptotically cylindrical $G_2$-manifold
asymptotic to $X\t(R,\iy)$, and $C$ an asymptotically cylindrical
coassociative 4-fold in $M$ asymptotic to $L\t(R',\iy)$.

\vspace{.1in}

 Let $\nu_L$ be the normal bundle of $L$ in $X$ and
$\exp_L:\nu_L\ra X$ be the exponential map. For $r>0$,
$B_r(\nu_L)$ is the subbundle of $\nu_L$ with fibre at $x$ the
open ball about 0 in $\nu_L\vert_x$ with radius $r$. Then for
small $\ep>0$, there is a tubular neighbourhood $T_L$ of $L$ in
$X$ such that $\exp_L:B_{2\ep}(\nu_L)\ra T_L$ is a diffeomorphism.
Also, $\nu_L\t\R\ra L\t\R$ is the normal bundle to $L\t\R$ in
$X\t\R$ with exponential map $\exp_L\t{\mathop{\rm
id}}:\nu_L\t\R\ra X\t\R$. Then $T_L\t\R$ is a tubular neighborhood
of $L\t\R$ in $X\t\R$, and $\exp_L\t{\mathop{\rm id}}:B_{2\ep}
(\nu_L)\t\R\ra T_L\t\R$ is a diffeomorphism.

\vspace{.1in}

Let $K,R$, $\Psi:X\t(R,\iy)\ra M\sm K$, and $K',R'>R$,
$\Phi:L\t(R',\iy)\ra C\sm K'$, and the normal vector field $v$ on
$L\t(R',\iy)$ be as before so that Diagram \eq{co2eq3} in Section
\ref{subsec2} commutes. Then $v$ is a section of
$\nu_L\t(R',\iy)\ra L\t(R',\iy)$, decaying at rate $O(e^{\be t})$
and by making $K'$ and $R'$ larger if necessary, we can suppose
the graph of $v$ lies in $B_\ep(\nu_L)\t(R',\iy)$.

\vspace{.1in}

Let $\pi:B_\ep(\nu_L)\t(R',\iy)\ra L\t(R',\iy)$ be the natural
projection. Then we define a diffeomorphism
\begin{equation}
\Xi:B_\ep(\nu_L)\t(R',\iy)\ra M \quad\text{by}\quad \Xi:w\mapsto
\Psi\bigl[(\exp_L\t{\mathop{\rm id}}) (v\vert_{\pi(w)}+w)\bigr].
\label{co4eq2}
\end{equation}
\noindent where $w$ is a point in $B_\ep(\nu_L)\t(R',\iy)$, in the
fibre over $\pi(w)\in L\t(R',\iy)$. Under the identification of
$L\t(R',\iy)$ with the zero section in $B_\ep(\nu_L)\t(R',\iy)$,
$\Xi\vert_{L\t(R',\iy)}\equiv\Phi$. Using $\Xi$ we can then define
an isomorphism $\xi$ between the vector bundles $\nu_L\t(R',\iy)$
and $\Phi^*(\nu_C)$ over $L\t(R',\iy)$, where $\nu_C$ is the
normal bundle of $C$ in $M$. This leads to the construction of
another diffeomorphism $\Th:B_{\ep'}(\nu_C)\ra T_C$ for
appropriate choices of a small $\ep'>0$, and a tubular
neighborhood $T_C$ of $C$ in $M$.

\vspace{.1in}

By choosing the local identification $\Th$ between $\nu_C$ and $M$
near $C$ that is compatible with the asymptotic identifications
$\Phi,\Psi$ of $C$, $M$ with $L\t\R$ and $X\t\R$ we can then
identify submanifolds $\ti C$ of $M$ close to $C$ with small
sections  of $\nu_C$. More importantly, the asymptotic convergence
of $\ti C$ to $C$, and so to $L\t\R$, is reflected in the
asymptotic convergence of sections of $\nu_C$ to 0.

\vspace{.1in}

We then define a map
$Q:L^p_{l+2,\ga}\bigl(B_{\ep'}(\La^2_+T^*C)\bigr)\ra \{$3-forms on
$C\}$ by $Q(\ze^2_+)=(\Th\ci\ze^2_+)^*(\vp)$ for $p>4$ and $l\ge
1$. That is, we regard the section $\ze^2_+$ as a map $C\ra
B_{\ep'} (\La^2_+T^*C)$, so $\Th\ci\ze^2_+$ is a map $C\ra
T_C\subset M$, and thus $(\Th\ci\ze^2_+)^*(\vp)$ is a 3-form on
$C$. Therefore if $\Ga_{\ze^2_+}$ is the graph of $\ze^2_+$ in
$B_{\ep'}(\La^2_+T^*C)$ and $\ti C=\Th(\Ga_{\ze^2_+})$ its image
in $M$, then $\ti C$ is coassociative if and only if
$\vp\vert_{\ti C}\equiv 0$, which holds if and only if
$Q(\ze^2_+)=0$. So basically $Q^{-1}(0)$ parametrizes
coassociative 4-folds $\ti C$ close to $C$. It turns out that
$Q:L^p_{l+2,\ga}\bigl(B_{\ep'}(\La^2_+T^*C) \bigr)\longra
L^p_{l+1,\ga}(\La^3T^*C)$ is a smooth map of Banach manifolds and
the linearization of $Q$ at $0$ is $\d
Q(0):\ze^2_+\mapsto\d\ze^2_+$.

\vspace{.1in}

As in the proof of McLean's Theorem, we then augment $Q$ by a
space of 4-forms on $C$ to make it elliptic and define
\begin{gather*}
P:L^p_{l+2,\ga}\bigl(B_{\ep'}(\La^2_+T^*C) \bigr)\t
L^p_{l+2,\ga}(\La^4T^*C)\longra
L^p_{l+1,\ga}(\La^3T^*C)\\
\text{by}\qquad P(\ze^2_+,\ze^4)=Q(\ze^2_+)+\d^*\ze^4
\end{gather*}

\noindent for $p>4$ and $l\ge 1$.

\vspace{.1in}

As mentioned before in Section \ref{subsec3}, on a noncompact
manifold, the ellipticity of a differential operator $A$ is not
sufficient to ensure that $A$ is Fredholm. Using the analytical
framework developed by Lockhart and McOwen in \cite{LM1} and
\cite{LM2}, we showed that in our case $\d P$ is not Fredholm if
and only if either $\ga=0$, or $\ga^2$ is a positive eigenvalue of
$\De=\d^*\d$ on functions on $L$, or $\ga$ is an eigenvalue of
$-*\d$ on coexact $1$-forms on $L$. Therefore we take $\ga$
sufficiently small to guarantee that $\d P$ is Fredholm.

\vspace{.1in}

We then show that $\Ker\bigl((\d_++\d^*)^p_{l+2,\ga}\bigr)$ is a
vector space of smooth, closed, self-dual $2$-forms and
$\Coker\bigl((\d^*_++\d)^q_{l+2,\ga}\bigr)$ is a vector space of
smooth, closed and coclosed $3$-forms. It turns out that the map
$\Ker\bigl((\d_++\d^*)^p_{l+2,\ga}\bigr)\ra H^2(C,\R)$,
$\chi\mapsto[\chi]$ induces an isomorphism of
$\Ker\bigl((\d_++\d^*)^p_{l+2,\ga}\bigr)$ with a maximal subspace
$V_+$ of the subspace $V\subseteq H^2(C,\R)$ on which the cup
product $\cup:V\t V\ra\R$ is positive definite. Hence
\begin{equation}
\dim\Ker\bigl((\d_++\d^*)^p_{l+2,\ga}\bigr)=\dim V_+,
\label{co3eq20}
\end{equation}
which is a topological invariant of $C,L$.

\vspace{.1in}

We finally show that $P$ maps
$L^p_{l+2,\ga}\bigl(B_{\ep'}(\La^2_+T^*C) \bigr)\t
L^p_{l+2,\ga}(\La^4T^*C)$ to the image of
$(\d_++\d^*)^p_{l+2,\ga}$ and the image of $Q$ consists of exact
3-forms. The Implicit Mapping Theorem for Banach spaces then
implies that $P^{-1}(0)$ is smooth, finite-dimensional and locally
isomorphic to $\Ker\bigl((\d_++\d^* )^p_{l+2,\ga}\bigr)$. As
$P^{-1}(0)=Q^{-1}(0)\t\{0\}$ we conclude that $Q^{-1}(0)$ is also
smooth, finite-dimensional and locally isomorphic to
$\Ker\bigl((\d_++\d^*)^p_{l+2,\ga}\bigr)$. Moreover, $Q^{-1}(0)$
is independent of $l$, and so consists of smooth solutions. This
proves Theorem \ref{co1thm}.

\end{proof}

\section{Coassociative Deformations with Moving Boundary and Proof of Theorem \ref{co1thm1}.}
\label{yco3}

We now prove Theorem \ref{co1thm1}. Let $(M,\vp,g)$ be an
asymptotically cylindrical $G_2$-manifold asymptotic to $X\t(R,\iy
)$, $R>0$, with decay rate $\al<0$. Let $C$ be an asymptotically
cylindrical coassociative 4-fold in $X$ asymptotic to $L\t(R',\iy
)$ for $R'>R$ with decay rate $\be$ for $\al\le\be<0$.

\vspace{.1in}

As in the proof of Theorem \ref{co1thm}, we suppose $\ga$
satisfies $\be<\ga<0$, and $(0,\ga^2]$ contains no eigenvalues of
the Laplacian $\De_L$ on functions on $L$, and $[\ga,0)$ contains
no eigenvalues of the operator $-*\d$ on coexact 1-forms on $L$.
The reason for this assumption is that, as in the fixed boundary
case, by Propositions \ref{comb1}, \ref{comb2}, we also have the
same linearized operator given as

\begin{equation} (\d_++\d^*)^p_{l+2,\ga}: F\times(
L^p_{l+2,\ga}(\La^2_+T^*C)\op L^p_{l+2,\ga}(\La^4T^*C)) \longra
L^p_{l+1,\ga}(\La^3T^*C) \label{co4eq1}
\end{equation}

\vspace{.1in}

\noindent for $p>4$ and $l\ge 1$ as in \S\ref{yco2}.

\vspace{.1in}

The conditions on $\ga$ imply that $[\ga,0)\cap\D_{(\d_+
+\d^*)_0}=\emptyset$ and hence $\ga\notin\D_{(\d_++\d^*)_0}$, so
that $(\d_++\d^*)^p_{l+2,\ga}$ is Fredholm. Here $\D_{(\d_+
+\d^*)_0}$ is the discrete set derived from the linearized
operator as in \cite{Salu0}.

\vspace{.1in}

Now, let $K,R$, $\Psi:X\t(R,\iy)\ra M\sm K$, and $K',R'>R$,
$\Phi:L\t(R',\iy)\ra C\sm K'$, and the normal vector field $v$ on
$L\t(R',\iy)$ be as before so that Diagram \eq{co2eq3} in Section
2 commutes. Let $\nu_C$ is the normal bundle of $C$ in $M$ and
$\nu_L$ be the normal bundle of $L$ in $X$. Then as in
\cite{Salu0}, by choosing an appropriate local identification
$\Th$ between $\nu_C$ and $M$ near $C$ that is compatible with the
asymptotic identifications $\Phi,\Psi$ of $C$, $M$ with $L\t\R$
and $X\t\R$ we can identify submanifolds $\ti C$ of $M$ close to
$C$ with small sections  of $\nu_C$.

\vspace{.1in}

Let $F\subset\R^d$ be an open subset of the space of special
Lagrangian deformations of $L$ in $(X,\omega, \Omega, g_X)$. Let
$L_0\in F$ be the starting point and $L_s$ be nearby special
Lagrangian submanifolds for some $s$ close to 0 in $F$. By McLean,
\cite{McLe}, the space of deformations $\mathcal M_{L_0}$ of a
special Lagrangian submanifold $L_0$ can be identified with closed
and coclosed 1-forms on $L_0$ and so $\mathcal M_{L_0}$
corresponds to $H^1(L_0,\R)$.

\vspace{.1in}

Moreover, one can naturally parametrize $L_s$ by $s\in
H^2(L_0,\R)$ and relate $\mathcal M_{L_0}$ to $H^2(L_0,\R)$. Given
a compact special Lagrangian $L_0$, let $U$ be a connected and
simply connected open neighbourhood of $L_0$ in $\mathcal
M_{L_0}$. One can construct natural local diffeomorphisms
$\mathcal A:\mathcal M_{L_0}\rightarrow H^1(L_0,\R)$ and $\mathcal
B:\mathcal M_{L_0}\rightarrow H^2(L_0,\R)$ as follows:

\vspace{.1in}

Let $L_s\in U$ be a special Lagrangian submanifold with phase $i$.
Then there exists a smooth path $\tilde\gamma:[0,1]\rightarrow U$
with $\tilde\gamma(0)=L_0$, and $\tilde\gamma(1)=L_s$ which is
unique up to isotopy. $\tilde\gamma$ parametrizes a family of
submanifolds of $X$ diffeomorphic to $L_0$, and can be lifted to a
smooth map $\Gamma:L_0\times [0,1]\rightarrow X$ with $\Gamma
(L_0\times\{s\})=\tilde\gamma(s)$. Now let $\Gamma^*(\omega)$ be a
2-form on $L_0\times[0,1]$. Then $\Gamma^*(\omega)|_{L_0\times
\{s\}}\equiv 0$ for each $s\in [0,1]$ as each fiber $\tilde
\gamma(s)$ is Lagrangian. Hence $\Gamma^*(\omega)$ can be written
as $\Gamma^*(\omega)=\alpha_s\wedge ds$ where $\alpha_s$ is a
closed 1-form on $L_0$ for $s\in[0,1]$. Then we integrate
$\alpha_s$ with respect to $s$ and take the cohomology class of
the closed 1-form $\int_0^1\alpha_s ds$. So we can define
$\mathcal A(L_s)=[\int_0^1\alpha_s ds]\in H^1(L_0,\R)$. Similarly,
we can write $\Gamma^*({\rm Im}(\Omega))=\beta_s\wedge ds$ where
$\beta_s$ is a closed 2-form on $L_0$ for $s\in[0,1]$ and define
$\mathcal B(L_s)=[\int_0^1\beta_s ds]\in H^2(L_0,\R)$. For more on
the constructions of the $H^2(L_0,\R)$ coordinates on the moduli
space of special Lagrangian submanifolds see \cite{GHJ} and
\cite{Hitc}.

\vspace{.1in}

We choose a smooth family of diffeomorphisms
$\vartheta_s:L_s\rightarrow L_0$ such that $\vartheta_0={\rm
id}_{L_0}$ and $L_s\cong L_0$ for small $s$. Identify a tubular
neighbourhood of $L_0$ in $X$ with a neighbourhood of the
zero-section in $T^*L_0$, then $L_s$ is the graph of $\zeta_s$,
where $\zeta_s$ is a small closed and coclosed 2-form on $L_0$
such that $[\zeta_s]=s$ in $H^2(L_0,\R)$.

\vspace{.1in}

Having got $\zeta_s$, for some small $s\in F$ we choose sections
$\varr_s$ of $\La^2_+T^*C$, that depend smoothly on $s$ where
$\varr_0=0$ and such that $\varr_s$ is asymptotic to
$\vars_s=\zeta_s+dt\wedge*_{L_0}(\zeta_s)$ where $*_{L_0}$ is the
star operator in $L_0$.

\vspace{.1in}

One way to do this is to take a smooth function $h:\R\rightarrow
[0,1]$ defined as $h(x)=0$ for $x\leq R$ and $h(x)=1$ for $x\geq
R+1$ and set $\varr_s=\vars_s\cdot h(x)$. Then $\varr_s=\vars_s$
in $(R+1,\infty)\t X$ and $\varr_s=0$ in $K$, where
$M=K\amalg(R,\infty)\t X$.

\vspace{.1in}

Given $\zeta_s$ as a section of $T^*L_0$, we can identify $T^*L_0$
with $\La^2_+(L_0\t\R)$ and identify this with $\La^2_+(C_0)$ near
$\infty$. Then $\tilde C_s$ is the graph of $\varr_s$, which is a
4-fold in $M$ asymptotic to $L_s\t\R$ but not necessarily
coassociative. So we obtain a family of 4-folds $\tilde C_s$ such
that $\tilde C_0=C_0=C$, which is a coassociative submanifold of
$M$, and $\tilde C_s$ is asymptotic to $L_s\t\R$, where $L_s$ is
special Lagrangian submanifold of $X$ as in Figure 1.

\vspace{.1in}

\begin{figure}[h]
\begin{center}

\includegraphics{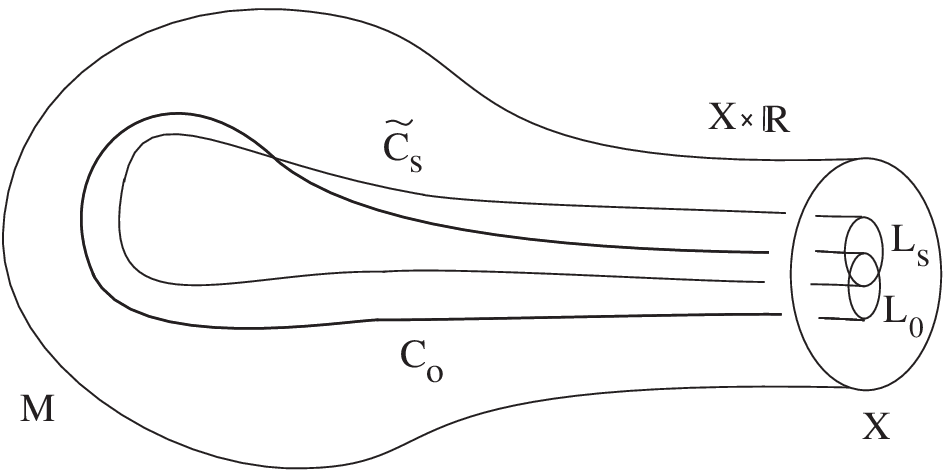}
\caption{ } \label{fig1}
\end{center}
\end{figure}

Note that when $\tilde C_s$ is asymptotic to $L_s\t(R,\infty)$
near $\infty$ then $\vp|_{\tilde C_s}\equiv O(e^{\ga t})$ near
$\infty$ (see Proposition \ref{co4prop1}). On the other hand,
since $L_s$ is a special Lagrangian, we can take $\varphi=0$ in a
small neighbourhood of the boundary $L_s$ in $L_s\t(R,\infty)$ and
so $\vp|_{\tilde C_s}$ is compactly supported and can take
$[\vp|_{\tilde C_s}]\in H^3_{\rm cs}(\tilde C_s,\R)$. Moreover,
since $L_0$ and $L_s$ are isotopic (which follows from our choice
of coordinates in $H^2(L_0,\R)$) we can take $[\vp|_{\tilde
C_s}]\in H^3_{\rm cs}(C_0,\R)$ which is independent of the choice
of $\tilde C_s$. And $[\vp|_{\tilde C_s}]$ should be zero for
there to exist a coassociative 4-fold $\tilde C_s$.

\vspace{.1in}

Now define $Q:F\times
L^p_{l+2,\ga}\bigl(B_{\ep'}(\La^2_+T^*C)\bigr)\ra \{$3-forms on
$C\}$ by $Q(s,\ze^2_+)=(\Th\ci(\ze^2_++\varr_s))^*(\vp)$. That is,
we regard the sum of sections $\ze^2_++\varr_s$ as a map $C\ra
B_{\ep'} (\La^2_+T^*C)$, so $\Th\ci(\ze^2_++\varr_s)$ is a map
$C\ra T_C\subset M$, and thus $(\Th\ci(\ze^2_++\varr_s))^*(\vp)$
is a 3-form on $C$. The point of this definition is that if
$\Ga_{(\ze^2_++\varr_s))}$ is the graph of $\ze^2_++\varr_s$ in
$B_{\ep'}(\La^2_+T^*C)$ and $\ti C_s=\Th(\Ga_{(\ze^2_++\varr_s)})$
its image in $M$, then $\ti C_s$ is coassociative if and only if
$\vp\vert_{\ti C_s}\equiv 0$, which holds if and only if
$Q(s,\ze^2_+)=0$. So $Q^{-1}(0)$ parametrizes coassociative
4-folds $\ti C_s$ close to $C$.

\vspace{.1in}

Note that to form $\Gamma(\ze^2_++\varr_s)$ we have chosen a
diffeomorphism $U\subseteq \La^2_+T^*C_0 =\nu_{C_0} \rightarrow
V\subseteq M$ where $U$ and $V$ are asymptotic to $U_0\subset
\La^2_+T^*(L_0\t\R)$ and $V_0\subset X\t\R$, respectively. Also,
we need $s$ to be sufficiently small in $F$ so that $L_s\times \R$
lies in $V_0$, in other words $L_s\t\R$ is the graph of $\vars_s$
in $\La^2_+T^*(L_0\t\R)$.

\vspace{.1in}

Hence we consider $\ze^2_+ +\varr_s$ for $\ze^2_+ \in
L^p_{l+2,\ga}(\La^2_+T^*C)$ and $\varr_s\in
C^{\infty}(\La^2_+T^*C)$ and define
\begin{equation*}
\begin{gathered}
Q:F\t L^p_{l+2,\ga}\bigl(\La^2_+T^*C )\ra L^p_{l+1,\ga}(\La^3T^*C),\\
\text{by}\qquad Q(s,\ze^2_+)=\pi_*(\vp\vert_{\Ga(\ze^2_+
+\varr_s)}).
\end{gathered}
\end{equation*}

\begin{prop} $Q:F\times L^p_{l+2,\ga}\bigl(B_{\ep'}(\La^2_+T^*C)
\bigr)\longra L^p_{l+1,\ga}(\La^3T^*C)$ is a smooth map of Banach
manifolds. The linearization of $Q$ at $0$ is $\d
Q(0,0):(s,\ze^2_+)\mapsto\d(\ze^2_++\varr_s )$. \label{co4prop1}
\end{prop}

\begin{proof}
The functional form of $Q$ is
\begin{equation*}
Q(s,\ze^2_+)\vert_x=H\bigl(s,
x,\ze^2_+\vert_x,\na\ze^2_+\vert_x\bigr) \quad\text{for $x\in C$,}
\end{equation*}

\noindent where $H$ is a smooth function of its arguments. Since
$p>4$ and $l\ge 1$ by Sobolev embedding theorem we have
$L^p_{l+2,\ga}(\La^2_+T^*C)\hookrightarrow C^1_\ga(\La^2_+T^*C)$.
General arguments then show that locally $Q(s,\ze^2_+)$ is
$L^p_{l+1}$.

\vspace{.1in}

From \cite{Salu0}, we know that $Q(0,\ze^2_+)$ lies in
$L^p_{l+1,\ga}(\La^3T^*C)$. When $s\neq0$, then
$Q(s,0)=\pi_*(\varphi|_{\Gamma(\varr_s)})$. By construction
$\Gamma(\varr_s)$ is asymptotic to $L_s\t\R$ which is
coassociative in $X\t\R$ as $L_s$ is a special Lagrangian 3-fold.
Identify $M\setminus K$ with $X\t\R$. Then $\varphi$ on
$M\setminus K$ can be written as $\varphi=\omega\wedge dt
+Re(\Omega)+O(e^\alpha t)$ where $\alpha$ is the rate for $M$
converging to $X\t\R$. $\Gamma(\varr_s)$ is the graph of $\varr_s$
which is equal to $L_s\t (R+1,\infty)$ in $X\t(R+1,\infty)$. Since
$L_s$ is a special Lagrangian submanifold with phase $i$, we get
for $t>R+1$, $\varphi|_{\Gamma(\varr_s)}=(\omega\wedge dt
+Re(\Omega)+O(e^\alpha t))|_{L_s\t(R+1,\infty)}=0+O(e^\alpha
t)|_{L_s\t(R+1,\infty)}$.

\vspace{.1in}

Therefore for $\Gamma(\varr_s)$ the error term
$\varphi|_{\Gamma(\varr_s)}$ comes from the degree of the
asymptotic decay. In particular, as $\alpha<\gamma$ we can assume
$\varphi|_{\Gamma(\varr_s)}\equiv O(e^{\gamma t})$. We can easily
arrange to choose $\varr_s$ such that
$\pi_*(\varphi|_{\Gamma(\varr_s)})\in L^p_{l+1,\ga}(\La^3T^*C)$
and that
$||\pi_*(\varphi|_{\Gamma(\varr_s)})||_{L^p_{l+1,\ga}}\leq c |s|$
for some constant $c$. This implies that $Q(s,\ze^2_+)$ lies in
$L^p_{l+1,\ga}(\La^3T^*C)$.

\vspace{.1in}

Finally, since we work locally, the linearization of $Q$ is again
$\d$ as before.
\end{proof}

Next we show that the image of $Q$ consists of exact 3-forms in
$L^p_{l,\ga}(\La^3T^*C)$.

\begin{prop}
$Q\bigl(F\t L^p_{l+2,\ga}(B_{\ep'}(\La^2_+T^*C))
\bigr)\!\subseteq\!\d\bigl(L^p_{l+1,\ga}(\La^2T^*C)\bigr)
\!\subset\!L^p_{l,\ga}(\La^3T^*C)$. \label{co14prop2}
\end{prop}

\begin{proof}

This will be an adaptation of a similar proof in \cite{Salu0}.
Consider the restriction of the 3-form $\vp$ to the tubular
neighborhood $T_C$ of $C$. As $\vp$ is closed, and $T_C$ retracts
onto $C$, and $\vp\vert_C\equiv 0$, we see that $\vp\vert_{T_C}$
is exact. Thus we may write $\vp\vert_{T_C}= \d\th$ for $\th\in
C^\iy(\La^2T^*T_C)$. Since $\vp\vert_C\equiv 0$ we may choose
$\th\vert_C\equiv 0$. Also, as $\vp$ is asymptotic to $O(e^{\be
t})$ with all its derivatives to a translation-invariant 3-form
$\vp_0$ on $X\t\R$, we may take $\th$ to be asymptotic to
$O(e^{\be t})$ with all its derivatives to a translation-invariant
2-form on~$T_L\t\R$. By Proposition \ref{co4prop1}, the map
$(s,\ze^2_+)\mapsto(\Th\ci(\ze^2_+ +\varr_s))^*(\th)$ maps
$F\times L^p_{l+2,\ga}\bigl(B_{\ep'}(\La^2_+T^*C)\bigr)\ra
L^p_{l+1,\ga}(\La^3T^*C)$. But

\begin{equation*}
Q(s,\ze^2_+)=(\Th\ci(\ze^2_+ +\varr_s))^*(\vp)=(\Th\ci(\ze^2_+
+\varr_s))^*(\d\th)=
\d\bigl[(\Th\ci(\ze^2_++\varr_s))^*(\th)\bigr],
\end{equation*}

\noindent so we can conclude that
$Q(s,\ze^2_+)\in\d\bigl(L^p_{l+1,\ga}(\La^2T^*C)\bigr)$ for
$\ze^2_+\in L^p_{l+2,\ga}(B_{\ep'}(\La^2_+T^*C))$ and $\varr_s\in
C^{\infty}(\La^2_+T^*C)$.
\end{proof}

As in \cite{Salu0}, we now augment $Q$ by a space of 4-forms on
$C$ to make it elliptic. That is, we define
\begin{equation*}
\begin{gathered}
P:F\t L^p_{l+2,\ga}\bigl(\La^2_+T^*C \oplus \La^4T^*C\bigr)\ra L^p_{l+1,\ga}(\La^3T^*C),\\
\text{by}\qquad P(s,\ze^2_+,\ze^4)=\pi_*(\vp\vert_{\Ga(\ze^2_+
+\varr_s)})+\d^*\ze^4.
\end{gathered}
\end{equation*}

Note that by the same discussion as in Proposition \ref{co4prop1},
we can also take that $P(s,\ze^2_+,\ze^4)$ lies in
$L^p_{l+1,\ga}(\La^3T^*C)$.

\vspace{.1in}

Proposition \ref{co4prop1} implies that the linearization $\d P$
of $P$ at 0 is the Fredholm operator $(\d_++\d^*)^p_{l+2,\ga}$ of
\eq{co4eq1}. Define $\mathcal C$ to be the image of
$(\d_++\d^*)^p_{l+2,\ga}$. Then $\mathcal C$ is a Banach subspace
of $L^p_{l+1,\ga}(\La^3T^*C)$, since $(\d_++\d^*)^p_{l+2,\ga}$ is
Fredholm. We show $P$ maps into $\mathcal C$.

\begin{prop} $P$ maps $F\t L^p_{l+2,\ga}\bigl(\La^2_+T^*C \oplus \La^4T^*C\bigr)\longra\mathcal C$.
\label{co4prop3}
\end{prop}

\begin{proof} Let $s\in F$ and $(\ze^2_+,\ze^4)\in L^p_{l+2,\ga}
\bigl(B_{\ep'}(\La^2_+T^*C)\bigr)\t L^p_{l+2,\ga} (\La^4T^*C)$, so
that $P(s,\ze^2_+,\ze^4)=Q(s,\ze^2_+)+\d^*\ze^4$ lies in
$L^p_{l+1,\ga}(\La^3T^*C)$. We must show it lies in $\mathcal C$.
Since $\ga\notin\D_{(\d_++\d^*)_0}$, this holds if and only if
\begin{equation}
\big\langle Q(s,\ze^2_+)+\d^*\ze^4,\chi\big\rangle_{L^2}=0
\qquad\text{for all $\chi\in\Ker\bigl((\d^*_++\d)^q_{
m+2,-\ga}\bigr)$,} \label{co4eq6}
\end{equation}

\noindent where $\frac{1}{p}+\frac{1}{q}=1$ and~$m\ge 0$.


\vspace{.1in}

We proved earlier that $\Ker\bigl((\d^*_++\d)^q_{ m+2,-\ga}\bigr)$
consists of closed and coclosed 3-forms $\chi$. We also know that
$Q(s,\ze^2_+)=\d\la$ for $\la\in L^p_{l+1,\ga}(\La^2T^*C)$ by
Proposition \ref{co14prop2}. So

\begin{equation*}
\big\langle Q(s,\ze^2_+)+\d^*\ze^4,\chi\big\rangle_{L^2}=
\big\langle \d\la,\chi\big\rangle_{L^2}+ \big\langle
\d^*\ze^4,\chi\big\rangle_{L^2}= \big\langle
\la,\d^*\chi\big\rangle_{L^2}+ \big\langle
\ze^4,\d\chi\big\rangle_{L^2}=0
\end{equation*}

\noindent for $\chi\in\Ker\bigl((\d^*_++\d)^q_{m+2,-\ga}\bigr)$,
as $\chi$ is closed and coclosed, and the inner products and
integration by parts are valid because of the matching of rates
$\ga,-\ga$ and $L^p,L^q$ with $\frac{1}{p}+ \frac{1}{q}=1$. So
\eq{co4eq6} holds, and $P$ maps into $\mathcal C$.
\end{proof}

\vspace{.1in}

Proposition \ref{co4prop3} implies that we can now apply Banach
Space Implicit Mapping Theorem, and therefore conclude that
$P^{-1}(0)$ is smooth, finite-dimensional and locally isomorphic
to $\mathcal{A}=\Ker\bigl((\d_++\d^* )^p_{l+2,\ga}\bigr)\subset
F\times L^p_{l+2,\ga}(\La^2_+T^*C \op \La^4T^*C)$.

\vspace{.1in}

Note that our original map was $Q$ and we needed the extra term
from the space of 4-forms on $C$ just to make $Q$ elliptic.
Therefore, we need to show that the following lemma holds.

\vspace{.1in}

\begin{lem} $P^{-1}(0)=Q^{-1}(0)\t\{0\}$.
\label{co4lem1}
\end{lem}

\begin{proof}  Let's assume that  $(s,\ze^2_+,\ze^4)\in P^{-1}(0)$, so that $Q(s,\ze^2_+)
+\d^*\ze^4=0$. This should imply $\ze^4=0$, so that
$Q(s,\ze^2_+)=0$, and therefore $P^{-1}(0)\subseteq
Q^{-1}(0)\t\{0\}$. By Proposition \ref{co14prop2} we have
$Q(s,\ze^2_+)=\d\la$ for $\la\in L^p_{l+1,\ga}(\La^2T^*C)$, so
$\d\la=-\d^*\ze^4$. Hence

\begin{equation*}
\lnm{\d^*\ze^4}{2}^2= \big\langle
\d^*\ze^4,\d^*\ze^4\big\rangle_{L^2}= -\big\langle
\d^*\ze^4,\d\la\big\rangle_{L^2}=
-\big\langle\ze^4,\d^2\la\big\rangle_{L^2}=0,
\end{equation*}

\noindent where the inner products and integration by parts are
valid as $L^p_{l+2,\ga}\hookrightarrow L^2_2$. Thus
$Q(s,\ze^2_+)=\d^*\ze^4=0$. But $\d^*\ze^4\cong\na\ze^4$ as
$\ze^4$ is a 4-form, so $\ze^4$ is constant. Since also $\ze^4\ra
0$ near infinity in $C$, we have $\ze^4\equiv 0$.

\vspace{.1in}

By construction, it is straightforward to show that
$Q^{-1}(0)\t\{0\}\subseteq P^{-1}(0)$.

\end{proof}

\begin{prop} Let $F\subset H^2(L,\R)$ be the subspace of special Lagrangian
deformations of the boundary $L$. Also let $\d
P_{(0,0,0)}(s,\ze^2_+,\ze^4)$ represent the linearization of the
deformation map $P$ at 0 with moving boundary and $\d
P^f_{(0,0)}(\ze^2_+,\ze^4)$ represent the linearization of the
deformation map $P^f$ at 0 with fixed boundary. Then
\begin{equation}
{\rm Index\; of\;} \d P_{(0,0,0)}= {\rm dim\;} F + {\rm index\;
of\;} \d P^f_{(0,0)}.
\end{equation}
\label{comb1}
\end{prop}

\begin{proof}
At $s=0$

\begin{equation}
P(0,\ze^2_+,\ze^4)=P^f(\ze^2_+,\ze^4)=\pi_*(\vp|_{\Gamma(\ze^2_+)})+\d
^*\ze^4
\end{equation}

Then at $s=0$, the linearization at $(0,0,0)$ is,

\begin{equation}
\d P_{(0,0,0)}(0,\ze^2_+,\ze^4)=\d P^f_{(0,0)}(\ze^2_+,\ze^4)
\end{equation}

\noindent and since $\d P_{(0,0,0)}$ is linear we have

\begin{equation}
\d P_{(0,0,0)}(s,\ze^2_+,\ze^4)=\d P_{(0,0,0)}(s,0,0)+\d
P_{(0,0,0)}(0,\ze^2_+,\ze^4)
\end{equation}

\noindent where $\d P_{(0,0,0)}(s,0,0)$ is finite dimensional,
$s\in T_0 F=\R^d$.

Then
\begin{equation}
{\rm Index\; of\;} \d P_{(0,0,0)}: F\t
L^p_{l+2,\ga}\bigl(\La^2_+T^*C \oplus \La^4T^*C\bigr)\ra
L^p_{l+1,\ga}(\La^3T^*C)
\end{equation}

\begin{equation}
 =  {\rm index\;of\;}
\d P^f_{(0,0)}: F\t L^p_{l+2,\ga}\bigl(\La^2_+T^*C \oplus
\La^4T^*C\bigr)\ra L^p_{l+1,\ga}(\La^3T^*C)
\end{equation}

\begin{equation}
= {\rm dim\;} F + ({\rm index\; of\;} \d P^f_{(0,0)}:
L^p_{l+2,\ga}\bigl(\La^2_+T^*C \oplus \La^4T^*C\bigr)\ra
L^p_{l+1,\ga}(\La^3T^*C)).
\end{equation}

\end{proof}

In Proposition \ref{comb2} we determine the set $F$. First, we
construct a linear map $\Up:H^2(L,\R)\rightarrow H^3_{\rm cs}(C)$
and then show that $F$ should be restricted to the kernel of
$\Up$.

\vspace{.1in}

\begin{prop} There is a linear map $\Up:H^2(L,\R)\rightarrow H^3_{\rm cs}(C)$ given explicitly as $ \Up(\vars):[ \vars ]\rightarrowtail [\d\tilde
\vars]$ and $F$ is an open subset of $\ker (\Up)$. \label{comb2}
\end{prop}

\begin{proof}

$F$ is a subset of the special Lagrangian 3-fold deformations of
$L$. So we could take $F\subset H^2(L,\R)$. We can construct a
linear map $\Up: H^2(L,\R)\rightarrow H^3_{\rm cs}(C)$ which comes
from the standard long exact sequence in cohomology:

\begin{equation}
\begin{gathered}
\xymatrix@C14pt@R17pt{ 0 \ar[r] & H^0(C) \ar[r] & H^0(L) \ar[r]
&H^1_{\rm cs}(C) \ar[r] & H^1(C) \ar[r] & H^1(L) \ar[r]
& H^2_{\rm cs}(C) \ar[d] \\
0 & H^4_{\rm cs}(C) \ar[l] & H^3(L) \ar[l] & H^3(C) \ar[l] &
H^3_{\rm cs}(C) \ar[l] & H^2(L) \ar[l] & H^2(C) \ar[l] }
\end{gathered}
\label{co3eq18}
\end{equation}

\noindent where $H^k(C)=H^k(C,\R)$ and $H^k(L)=H^k(L,\R)$ are the
de Rham cohomology groups, and $H^k_{\rm cs}(C,\R)$ is the
compactly-supported de Rham cohomology group.

\vspace{.1in}

The explicit definition of the map $\Up$ is that if $\vars$ is a
closed 2-form on $L$ then extend $\vars$ smoothly to a 2-form
$\tilde\vars$ on $C$ which is asymptotic to $\vars$ on
$L\t(\R,\infty)$. Then $\d \tilde\vars$ is a closed compactly
supported 3-form on $C$. So $\Up: H^2(L,\R)\rightarrow H^3_{\rm
cs}(C)$ can be defined explicitly as $\Up(\vars):[ \vars
]\rightarrowtail [\d\tilde \vars]$. Also note that by construction
of the coordinates in $H^2(L,\R)$, $[\vp|_{\tilde C_s}]=\Up([s])$,
for $s\in H^2(L,\R)$ and for the graphs of $\varr_s$.

\vspace{.1in}

 The potential problem here is that the image of $\Up$ may not
be zero in $H^3_{\rm cs}(C)$.  In \cite{Salu0}, we studied the
properties of the kernel and cokernel of the operator $d_++d^* $.
The dimension of the cokernel is

\begin{equation*}
\dim\Ker\bigl((\d^*_++\d)^q_{m+2,-\ga}\bigr)=b^0(L)-b^0(C)+b^1(C).
\end{equation*}

The map $\Ker\bigl((\d^*_++\d)^q_{m+2,-\ga}\bigr)\ra H^3(C,\R)$,
$\chi\mapsto[\chi]$ is surjective, with kernel of dimension
$b^0(L)-b^0(C)+b^1(C)-b^3(C)\ge 0$, which is the dimension of the
kernel of $H^3_{\rm cs}(C,\R)\ra H^3(C,\R)$. So we can think of
$\Ker\bigl((\d^*_++\d)^q_{m+2,-\ga}\bigr)$ as a space of closed
and coclosed 3-forms filling out all of $H^3_{\rm cs} (C,\R)$ and
$H^3(C,\R)$. In other words cokernel has a piece looking like the
kernel of $H^3_{\rm cs}(C,\R)\ra H^3(C,\R)$. This is equivalent to
the image of $\Up: H^2(L,\R)\rightarrow H^3_{\rm cs}(C)$ by
exactness of the long exact sequence \eq{co3eq18}.

\vspace{.1in}

So the problem here is that if we allow the boundary to move
arbitrarily then the image of $\Up$ may not be zero in $H^3_{\rm
cs}(C)$. This means that, for $s\in H^2(L,\R)$, the construction
of asymptotically cylindrical coassociative 4-fold asymptotic to
$L_s\t(R,\infty)$ is obstructed if and only if $\Up(s)\neq 0$.
Therefore the set $F$ is restricted to $\ker\Up\subset H^2(L,\R)$.

\end{proof}

\vspace{.1in}

Next, we determine the dimension of the moduli space $\M_C^{\ga}$.

\begin{prop} The dimension of $\M_C^\ga$, the moduli space of
coassociative deformations of an asymptotically cylindrical
coassociative submanifold $C$ asymptotic to $L_f\t (R,\infty)$,
$f\in F$, with decay rate $\ga$ is

\begin{equation}
{\rm dim} (\M_C^\ga)={\rm dim} (V_+) +
b^2(L)-b^0(L)+b^0(C)-b^1(C)+b^3(C).
\end{equation}
\end{prop}

\begin{proof}
Let $b^k(C)$, $b^k(L)$ and $b^k_{\rm cs}(C)$ be the corresponding
Betti numbers as before. The actual dimension of the moduli space
for free boundary is the sum of the dimension of the fixed
boundary and the dimension of the kernel of $\Up$.

We know from Propositions \ref{comb1}, and \ref{comb2} that
\begin{equation}
{\rm index\; of\;} \d P_{(0,0,0)}= {\rm dim\;} F + {\rm index\;
of\;} \d P^f_{(0,0)}.
\end{equation}

In particular,

\begin{equation}
{\rm dim}\; (\ker \d P_{(0,0,0)})= {\rm dim}\;(\ker \Up)+ {\rm
dim}\; (V_+).
\end{equation}

Taking alternating sums of dimensions in \eq{co3eq18}, shows that
the dimension of the kernel of $\Up$ is

\begin{equation}
\begin{split}
{\rm dim}\;(\ker\Up) &=b^2(L)- {\rm Im}(\Up)\\
&=b^2(L)-[b^0(L)-b^0(C)+b^1(C)-b^3(C)].
\end{split}
\end{equation}

Therefore the dimension of $\M_C^\ga$, the moduli space of
coassociative deformations of an asymptotically cylindrical
coassociative submanifold $C$ asymptotic to $L_s\times
(R,\infty)$, $s\in F$, with decay rate $\ga$ is

\begin{equation}
\begin{split}
{\rm dim} (\M_C^\ga)&={\rm dim} (V_+) + {\rm dim} (\ker\Up).\\
&={\rm dim} (V_+) + b^2(L)-b^0(L)+b^0(C)-b^1(C)+b^3(C).
\end{split}
\end{equation}

\end{proof}

\vspace{.1in}

And finally straightforward modifications of Propositions 4.6-4.8
in \cite{Salu0} provides us the standard elliptic regularity
results and the bootstrapping argument to conclude that
$Q^{-1}(0)$ is a smooth, finite dimensional manifold, and so we
complete the proof of Theorem \ref{co1thm1}. Here we will skip
these details to avoid repetition, for more on regularity results
see \cite{Salu0}.

\section{Applications: Topological Quantum Field Theory of Coassociative Cycles}
\label{yco4}

In \cite{Leun2}, it was mentioned that if one could show the
analytical details of the deformation theory of asymptotically
cylindrical coassociative submanifolds inside a $G_2$-manifold
then it would be possible to study global properties of these
moduli spaces. In particular, one needs a result like Theorem
\ref{co1thm1} which shows that $H^2_+(C,L)$ parametrizes the
coassociative deformations of $C$ with boundary $L$. In this
section using Theorem \ref{co1thm1} and a well-known theory of
(anti) self-dual connections (for more on the subject see
\cite{DK}, \cite{FU}) we will verify this claim.

\vspace{.1in}

Here are some basic definitions:

\vspace{.1in}

Let $(M,\vp,g)$ be an asymptotically cylindrical $G_2$-manifold
asymptotic to $X\t(R,\iy)$, and $C$ an asymptotically cylindrical
coassociative 4-fold in $M$ asymptotic to $L\t(R',\iy)$. Let $E$
be a fixed rank one vector bundle over $C$. A connection $D_E$ on
$C$ has finite energy if

$$\displaystyle \int_C |F_E|^2\;dV < \infty$$

\noindent where $F_E$ is the curvature of the connection $D_E$ and
$dV$ is the volume form with respect to the metric. Note that
connections with finite energy play an important role in
Yang-Mills Theory.

\begin{dfn}
Let $(X,\om, J, g_X, \Om)$ be a Calabi--Yau 3-fold and $L$ be a
special Lagrangian submanifold of $X$. Also let $(M,\vp,g_M)$ be a
$G_2$-manifold asymptotic to $X\t(R,\iy)$. A pair $(C,D_E)$ is
called a coassociative cycle if $C$ is a coassociative submanifold
of $M$ asymptotic to $L\t(R',\iy)$ and $D_E$ is an anti-self-dual,
unitary connection over $C$ with finite energy.
\end{dfn}

\begin{dfn}
Let $(X,\om,J,g_X,\Om)$ be a Calabi--Yau 3-fold and $L$ be a
special Lagrangian submanifold in $X$. A pair $(L,D_{E'})$ is
called a special Lagrangian cycle if $D_{E'}$ is a unitary flat
connection over $L$ induced from $D_E$.
\end{dfn}

Let $\mathcal{M}^{slag}(X)$ denote the moduli space of special
Lagrangian cycles in $X$. In \cite{Hitc}, Hitchin proved the
following theorem:

\begin{thm}
The tangent space to $\mathcal{M}^{slag}(X)$ is naturally
identified with the space $H^2(L,\R)\t H^1(L, ad(E'))$. For line
bundles over $L$, the cup product $\cup:H^2(L,\R)\t
H^1(L,\R)\longrightarrow \R$ induces a symplectic structure on
$\mathcal{M}^{slag}(X)$.
\end{thm}

Now using Theorem \ref{co1thm1} we verify the proof of the
following theorem which was claimed by Leung [Claim 10, Sec.4] in
\cite{Leun2}.

\begin{thm}
Let $X$ be a Calabi--Yau 3-fold and $L$ be a special Lagrangian
submanifold in $X$. Let $M$ be a $G_2$-manifold asymptotic to
$X\t(R,\iy)$, and $C$ a coassociative 4-fold in $M$ asymptotic to
$L\t(R',\iy)$. Let also $\mathcal{M}^{coas}(M)$ be the moduli
space of coassociative cycles in $M$ and $\mathcal{M}^{slag}(X)$
be the moduli space of special Lagrangian cycles in $X$. Then the
boundary map

\begin{equation}
\begin{split}
b:\mathcal{M}^{coas}(M)&\longrightarrow
\mathcal{M}^{slag}(X)\\
b(C, D_E)&=(L, D_{E'})
\end{split}
\end{equation}

\noindent is a Lagrangian immersion. \label{co1thm23}
\end{thm}

\begin{proof}

Let $\delta: H^1(L)\rightarrow H^2_+(C,L)$, $j^*: H^2_+(C,L)
\longrightarrow H^2_+(C)$, $i^*:H^2_+(C)\longrightarrow H^2(L)$ be
the canonical maps. They give dual maps $i_*:H_2(L)\longrightarrow
H^+_2(C)$, $j_*: H^+_2(C)\longrightarrow H^+_2(C,L)$ and $\partial
:H^+_2(C,L)\longrightarrow H_1(L)$. Then we have the following
long exact sequences:

\vspace{.1in}
\begin{equation*}
\begin{matrix}
\ldots \longrightarrow H^1(L)& \longrightarrow & H^2_+(C,L)& \longrightarrow & H^2_+(C)&\longrightarrow &H^2(L) &\longrightarrow &H^3(C,L) &\longrightarrow \ldots\\
\;\;\;\;\;\;\;\;\;\;\;\; \downarrow \cong &      & \downarrow \cong  &             & \downarrow \cong &            &\downarrow \cong&           & \downarrow \cong&                  \\
\ldots \longrightarrow H_2(L)& \longrightarrow  & H^+_2(C)
&\longrightarrow   & H^+_2(C,L)&\longrightarrow & H_1(L)&
\longrightarrow & H_1(C) &\longrightarrow  \ldots\\
\end{matrix}
\end{equation*}

\vspace{.1in}

Note that the classes coming from the boundary $L$ should have
zero self-intersection. So we can rewrite the long exact sequences
above starting from 0, instead of $H^1(L)$ and $H_2(L)$.

\vspace{.1in}

From the sequences we get
$$\ker\delta\cong\ker i_*\cong im \partial\cong H^+_2(C,L)/\ker\partial=H^+_2(C,L)/im(j_*)$$
\noindent which implies
\begin{equation}
\begin{split}
{\rm dim} \;(\ker\delta) &= {\rm dim}\; (H^+_2(C,L))-{\rm dim}(im (j_*))\\
\end{split}
\label{co2eq223}
\end{equation}

Note that the only classes that survive in
$H^+_2(C)\longrightarrow H^+_2(C,L)$, (i.e. do not go to zero)
have self intersection 0. So one can identify $im (j_*)$ with
$H^+_2(C)$.

\vspace{.1in}

This implies that there is a short exact sequence
$$0\longrightarrow im (j_*) \longrightarrow H^+_2(C,L)\longrightarrow im\partial \longrightarrow 0, $$

\noindent or equivalently

$$0\longrightarrow H^+_2(C) \longrightarrow H^+_2(C,L)\longrightarrow \ker\delta \longrightarrow 0.$$

\vspace{.1in}

Then as a consequence of \eq{co2eq223} and Theorem \ref{co1thm1}
we conclude that $H^2_+(C,L)=V_+\oplus \ker \delta$ parametrizes
the deformations of $C$ with moving boundary $\partial C$. From
our previous work in \cite{Salu0}, it follows that $H^2_+(C)\cong
V_+$ parametrizes the deformations of $C$ with fixed boundary and
McLean showed that $H^2(L)$ gives the special Lagrangian
deformations of $L$.

\vspace{.1in}

Note that the linearization of the boundary value map
$\mathcal{M}^{coas}(M)\longrightarrow \mathcal{M}^{slag}(X)$ in
Theorem \ref{co1thm23} is given by $i^*:H^2_+(C)\longrightarrow
H^2(L)$ and for the connection part $\beta: H^1(C)\longrightarrow
H^1(L)$. It is straightforward that Im$(i^*)\oplus$Im$\beta$ is a
subspace of $H^2(L)\oplus H^1(L)$. By definition of cup product,
the symplectic structure reduces to 0 on Im$(i^*)\oplus$Im$\beta$
and by Poincar$\acute{e}$ Duality,
dim(Im$(i^*)\oplus$Im$\beta)=\frac{1}{2}$ dim($H^2(L)\oplus
H^1(L)$).

\vspace{.1in}

Thus we conclude that Im$(i^*)\oplus$Im$\beta$ is a Lagrangian
subpace of $H^2(L)\oplus H^1(L)$ with the symplectic structure
defined above and conclude the proof of Theorem \ref{co1thm23}.

\end{proof}

{\small{\it Acknowledgements.} Special thanks to an anonymous
referee for many useful comments which greatly improved the
previous version of this paper. Also thanks to Dominic Joyce for
his help during the early stages of this project.}

\end{document}